\documentclass{amsart}

\newtheorem{theorem}{Theorem}[section]
\newtheorem{lemma}[theorem]{Lemma}
\newtheorem{proposition}[theorem]{Proposition}

\theoremstyle{definition}

\theoremstyle{remark}
\newtheorem{remark}[theorem]{Remark}

\numberwithin{equation}{section}

\def\de#1{\textit{#1}}
\def\td{{\tilde d}}
\def\tg{{\tilde g}}
\def\tJ{{\tilde J}}
\def\hg{{\hat g}}
\def\tx{x}
\def\ty{{y}}
\def\talpha{{\tilde \alpha}}
\def\ta{{\tilde \alpha}}
\def\circle{{S^1}}

\def\reals{{\mathbb R}}
\def\integers{{\mathbb Z}}
\def\naturals{{\mathbb N}}
\def\raw{\rightarrow}

\def\part{\Gamma}
\def\pKN{p_{K,N}}
\def\pKNN{p_{K+1,N}}
\def\pkn{p_{k,n}}

\def\hx{{\hat x}}
\def\hI{{\hat I}}
\def\hJ{{\hat J}}
\def\Ispec{[\pKN, \pKNN]}
\def\var{\text{var}}
\def\us{{\underline s}}
\def\ut{{\underline t}}
\def\htop{h_{top}}
\def\Inv{^{-1}}
\def\increasesto{\nearrow}

\begin{document}

\title{Semiconjugacies to angle-doubling}

\author{Philip Boyland}
\address{
Department of Mathematics,
University of Florida,
Gainesville, FL 32605-8105}
\email{boyland@math.ufl.edu}


\date{October 31, 2004.}


\begin{abstract}
A simple consequence of a
theorem of Franks says that whenever a continuous map, $g$, is 
homotopic to angle doubling on the circle it is semiconjugate to 
it. We show that
when this semiconjugacy has one disconnected point inverse,
then the typical point in the circle has a point inverse
with uncountably many connected components. Further, 
in this case the topological entropy of $g$ is strictly larger than
that of angle doubling, and the semiconjugacy  has 
unbounded variation.  An analogous theorem holds
for degree-$D$ circle maps with $D > 2$.
\end{abstract}

\maketitle

\section{Introduction}
The angle doubling map, $d$,  on the circle, $\circle := \reals/\integers$, 
is an oft cited example of a chaotic dynamical system.
If we define the itinerary of $\theta\in\circle$ 
as the sequence $\us$ defined by $s_i =  0$ if $0 < d^i(\theta) \leq 1/2$ and
$s_i =  1$ if $1/2 < d^i(\theta) \leq 1$, then for
{\it any} sequence of $0's$ and $1's$ we can find a $\theta$
which has that sequence as its itinerary. Thus the system embeds
the randomness of a sequence of coin tosses within its dynamics.

This dynamical complication of angle-doubling is actually
topological in character in the sense that
it cannot be removed by continuously deforming the system. 
A theorem of Franks (\cite{F}) shows that any circle map that is homotopic
to $d$ has dynamics at least as complicated as those of $d$ in
the precise sense given in the next theorem.
(Angle-doubling on a circle is a simplest case of a much more
general theorem.)

\begin{theorem}[Franks]\label{th1}
If $g$ is a continuous, circle map that is homotopic to the
angle-doubling map $d$, then there exists a continuous, onto map
 $\alpha:\circle\raw\circle$ with $\alpha\circ g = d\circ \alpha$.
\end{theorem}

An $\alpha$ as in the theorem is called a \de{semiconjugacy}
of $g$ to $d$. The theorem can be informally understood by noting
that whenever $g$ is homotopic to $d$, the map $g^n$ must of
necessity wrap the circle $2^n$ times around itself, and
so iterates of $g$ have an unavoidable topological
complication.

A useful interpretation of the theorem
considers the point inverses, $\alpha\Inv(\theta)$, as
``fibers'' over the points $\theta$. The dynamics of
$g$ can be then thought of as twisted product 
with the base point $\theta$ moved
according to $d$ while the fiber over $\theta$ is mapped
by $g$ to  the fiber over $d(\theta)$. Thus all the information
about how the dynamics of $g$ differ from those of
$d$ is contained in the nature of the point inverses of
$\alpha$ and in the way in which these point inverses are transformed into
each other by $g$. 

If $\alpha$ is homeomorphism, each $\alpha\Inv(\theta)$ is a
single point, and so $g$ and $d$ have the same dynamics.
The next simplest case is when each $\alpha\Inv(\theta)$ is
a connected set, and thus is a point or an interval. In this case the 
essential difference
between the dynamics of $g$ and $d$ is contained in the
dynamics on intervals, a much studied subject.
The case of interest here is when $\alpha$ has at least
one disconnected point inverse. In this case the dynamics of $g$ 
are much more complicated than those of $d$ in the
sense that the typical fiber, $\alpha\Inv(\theta)$, has
uncountably many connected components.

\begin{theorem}\label{th2}
If $g$ is a continuous circle map that is  homotopic to the
angle-doubling map $d$ and $\alpha$ is its 
semiconjugacy to $d$, 
then  the following are equivalent:
\begin{itemize}

\item[(a)] There exists a point $\theta\in\circle$ with
$\alpha\Inv(\theta)$ disconnected,

\item[(b)] There exists a full measure, dense, $G_\delta$-set $\Lambda\subset\circle$
so that $\theta\in \Lambda$ implies that $\alpha\Inv(\theta)$ has uncountably
many connected components,

\item[(c)] The map $\alpha$ is not of bounded variation.
\end{itemize}
\noindent Further, in this case the topological entropy of
$g$ is strictly larger than that of $d$, $\htop(g) > \htop(d) = \log(2)$.
\end{theorem}

Note that the existence of the semiconjugacy yields that 
$\htop(g) \geq \htop(d)$, so the content of the last statement
of the theorem is the strict inequality.  
From the point of view developed before the theorem, this 
conclusion indicates that
the action of $g$ in permuting the fibers $\alpha\Inv(\theta)$ 
has positive entropy.

We briefly remark on related work.
The case not included in the theorem, namely when the semiconjugacy
has connected point inverses, includes the situation where $g$ is a
covering map (see the first paragraph of the proof of Theorem~\ref{th7}). 
The semiconjugacies of degree-two covering maps 
have been much studied from an analytic point of view 
(see, for example, section II.2 in \cite{dMvS}, and the references therein).
Also, there is a theorem in symbolic dynamics concerning
a semiconjugacy between two transitive subshifts of finite
type which bears a resemblance to Theorems~\ref{th2} and~\ref{th7} 
(see Remark~\ref{r3}).
Finally, there are theorems analogous to Theorems~\ref{th2} and~\ref{th7}
which hold
for homeomorphisms of the two-torus 
which are isotopic to Anosov diffeomorphisms.
These will be the subject of a subsequent paper.

While we state and prove our results 
for degree-two maps, it will be clear that
virtually identical proofs yield the 
analogous theorems for degree-$D$ circle maps with $D > 2$.

\section{Preliminaries}
The circle $\circle$ has universal cover $\reals$,
and the phrases \de{lift} and \de{projection}
always mean lifts to and projections from this cover.
A circle map is said to have degree $D\in\integers$ if it
is homotopic to $\theta \mapsto D \theta$. In the special 
case of degree two, we write the  angle doubling map as
$d(\theta) = 2 \theta$, and it has lift $\td(x) = 2x$.
Note that a map $g:\circle\raw\circle$ has degree two if and only if 
its lift can be written as 
\begin{equation}\label{cc}
\tg = \td  + \varphi 
\end{equation}
with $\varphi(\tx + 1) = \varphi(\tx)$. 

Given a degree-two circle map $g$ with lift $\tg$, 
for each $D\in \integers$, let $B_D$ be the Banach space of all lifts of
continuous degree-$D$ circle maps with the sup norm, and 
define $F_D:B_D \raw B_D$ by $F(f) = (f \circ\tg)/2$. 
It is easy to see that $F_D$ is a contraction mapping whose
fixed point $\ta_D$ satisfies  $\ta_D\circ \tg = \td\circ \ta_D$, and 
so projecting to the circle for any $D\not = 0$, 
we obtain Theorem~\ref{th1}.  This proof shows that for each $D$
the semiconjugacy $\alpha_D$ is the unique 
continuous, degree-$D$ map which satisfies 
$\alpha_D \circ g = d\circ \alpha_D$. In this paper we will only consider
the case $D = 1$, and given a degree-two $g$ by 
\de{its semiconjugacy} we always mean $\alpha_1$, which will
henceforth be denoted $\alpha$.
If we begin the iteration of $F_1$ with the identity map, $id$, we obtain
\begin{equation}\label{aa}
\frac{\tg^n}{2^n} = F_1^n(id) \raw \ta 
\end{equation}
uniformly.

It is also useful to consider an operator that acts on
the periodic parts of the maps. If the given degree-two map
is as in \eqref{cc}
and $C$ is the Banach space of $1$-periodic functions
with the sup norm, then $G:C\raw C$ defined by
$G(\sigma) = (\varphi + \sigma\circ\tg)/2$ is
also a contraction mapping, and if its fixed point is
$\gamma$,  then the lift of the semiconjugacy is $\ta = id + \gamma$.
If we begin the iteration of
$G$ with the zero map ${\bf 0}$, we obtain that
\begin{equation}\label{bb}
G^n({\bf 0}) = \sum_{i=0}^{n-1} \frac{\varphi\circ\tg^i}{2^{i+1}}
\raw \gamma
\end{equation}
uniformly and so 
\begin{equation*}
\ta = id + \sum_{i=0}^{\infty} \frac{\varphi\circ\tg^i}{2^{i+1}},
\end{equation*}
as could have been confirmed directly.

The semiconjugacy gives  a uniform bound on the 
distance between the $\tg$-orbit of $x$ and the $\td$ orbit of
$\ta(x)$. 
Using the semiconjugacy and $\ta = id + \gamma$
\begin{equation}\label{gs}
|\tg^n(\tx) - 2^n \ta(\tx)| =
|\tg^n(\tx) - \ta(\tg^n(\tx))| \leq \|\gamma\| \leq \|\varphi\|,
\end{equation}
for all $n$, 
where for the last inequality we used \eqref{bb}. In the language of
\cite{H}, this says that the orbits $o(x, g)$ and $o(\alpha(x), d)$ globally
shadow, where for a given map $f$, the \de{orbit} of a point $x$ is
$o(x, f) := \{f^n{x} : n = 0, 1, \ \dots\}$. It is worth noting that
Theorem~\ref{th1} can also be proved by a slight alteration of
global shadowing proof of the semiconjugacies to pseudoAnosov maps 
given in \cite{H}.

Recall that a map is called \de{light} if every point preimage is 
totally disconnected and \de{monotone} if every point preimage
is connected. A theorem of Eilenberg and Whyburn (independently)
says that for any
continuous map $f: X \raw Y$ with $X$ and $Y$ compact metric
spaces, there exist a compact metric space $Z$, a 
continuous light map $\ell:Z\raw Y$ and a continuous
monotone map $m:X\raw Y$, so
that $f = \ell m$. The decomposition is particularly simple
in the case at hand, $X = Y = \circle$, for since
connected components of point inverses are always closed
intervals, $Z = \circle$, and the monotone map $m$ simply collapses
certain intervals to points.

To study semiconjugacies $\alpha$ with
disconnected point preimages it is useful at first to
ignore the monotone part of $\alpha$ and assume
that $\alpha$ is light. 
We shall see in the proof of Theorem~\ref{th2} that by collapsing
collections of invariant intervals, any
degree-two $g$ can be projected to a degree-two map whose
semiconjugacy is light. The next proposition gives various dynamical
characterizations of those $g$ whose  semiconjugacies are light maps.

Recall that a map $f$ on a space $X$ is \de{locally eventually onto}
(leo) if for any open set
$U$ there is an $n>0$ so that $f^n(U) = X$. A map is \de{transitive} if
it has a dense orbit. A well-known characterization of transitivity 
on compact metric spaces is
that for all open $U$ and $V$ there exists an $n>0$ so that $f^n(U) \cap
V \not = \emptyset$, and so clearly  leo implies transitivity. For
a one-dimensional system an interval $J$ is \de{periodic} if there 
exists an $n>0$ so that $f^n(J)\subset J$, and $J$ is \de{wandering} if
for all $i\not = j$, $i, j\geq 0$, $f^i(J) \cap f^j(J) = \emptyset$.
Here and throughout this paper the terminology \de{interval} always means
a compact, nontrivial interval.

\begin{proposition}\label{th3}
 If $g$ is a continuous degree-two circle map
the following are equivalent:
\begin{itemize}

\item[(a)] The semiconjugacy $\alpha$ of $g$ to $d$ is light,

\item[(b)] $g$ is locally eventually onto,

\item[(c)] $g$ is transitive,

\item[(d)] $g$ is light and has no periodic or wandering intervals.
\end{itemize}
\end{proposition}

\begin{proof} If $J$ is a nontrivial interval and $\alpha$
is light, then there must exist $x_1, x_2\in \tJ$ with $\ta(x_2) > \ta(x_1)$ 
and $\tJ$ a lift of $J$.
Thus we may find an $n>0$ with $2^n\ta(x_2) - 2^n\ta(x_1) > 1 + 
2 \|\varphi\|$ where $\tg(x)$ is   as in 
\eqref{cc}. Thus by \eqref{gs}, $\tg^n(x_2) - \tg^n(x_1) > 1$, and so
$g^n(J) = \circle$. Therefore, (a) implies (b), and as noted above
the theorem, (b) implies (c). Now assume  that $o(x, g)$ is dense.
If $\alpha$ was not light, then for some  nontrivial interval
$J$, $\alpha(J) = \theta_0$, a point. Since $o(x, g)$ is dense, 
there are $i\not = j$ with $g^i(x)\in J$ and $g^j(x) \in J$. 
Thus $\alpha(g^i(x)) = \alpha(g^j(x)) = \theta_0$ for $i \not = j$, 
and so by the semiconjugacy, $d^i(\alpha(x)) = d^j(\alpha(x))$, and
so $o(\alpha(x), d)$ is eventually periodic. On the other hand,
the continuity of the semiconjugacy implies that $o(\alpha(x), d)$ 
is dense since $o(x, g)$ is, a contradiction, and so (c) implies (a).

Now (b) clearly implies (d). We finish by proving the contrapositive
of (d) implies (a), so assume $\alpha$ is not light, and thus there is
some nontrivial interval $J$ with $\alpha(J) = \theta_0$.
Now if $g^n(J)$ is a point for some $n>0$ or if $J$ wanders, we are done.
So we are left with the case  when there is an $i > j$ with 
$g^i(J)\cap g^j(J)\not = \emptyset$. The semiconjugacy then
yields that $d^i(\theta_0) = \alpha(g^i(J)) = \alpha(g^j(J)) = d^j(\theta_0)$.
Thus if $\hJ$ is the connected component of $\alpha\Inv(d^j(\theta_0))$
which contains $g^j(J)$, we must have $g^{i-j}(\hJ) \subset \hJ$, and so 
$g$ has a periodic interval. \end{proof}

We shall make frequent use of standard results and techniques 
of one-dimensional dynamics without mention, but 
for the reader's convenience we state the following
fundamental lemma.  Recall that  $I$ \de{covers} $J$ 
means that $J\subset I$. For more information on one-dimensional
dynamics see \cite{ALM}, \cite{BC}, or \cite{dMvS}. The version of
the lemma we give essentially comes from \cite{MN}.

\begin{lemma}\label{th4}
Assume that $f:\reals\raw \reals$ is continuous.
\begin{itemize}

\item[(a)] If $f(J) $ covers $I$, then there is an interval $J'\subset J$
so that $f(J') = I$ and no interior point of $J'$ maps to the 
boundary of $I$ under $f$.

\item[(b)] If $\{J_i\}$ is a finite collection of intervals
such that $f(J_i)$ covers $J_{i+1}$ for all $i$, then there exists
an interval $J'\subset J_0$ with $f^i(J') \subset J_i$ for all $i$.
If $\{J_i\}$ is an countable collection, then 
there is a $y\in J_0$ with $f^i(y) \in J_i$ for all $i$.
\end{itemize}
\end{lemma}

\section{The main lemmas}

The first main lemma locates a copy of the
dynamics of $d$ inside the dynamics of $g$. It makes no
assumptions about the lightness or injectivity of
the semiconjugacy.

\vfill\eject

\begin{lemma}\label{th5}
Given a degree-two circle map $g$ with
semiconjugacy $\alpha$, 
for each $r\in\reals$ let $p_r = \min\{\ta\Inv(r)\}$.
\begin{itemize}

\item[(a)] If $\tx < p_r$, then $\talpha(\tx) < r$.

\item[(b)] The map $r \mapsto p_r$ is order preserving.
 
\item[(c)] Each $p_r$ satisfies  $\tg(p_r) = p_{2r}$.

\item[(d)]  If $\tx < p_r$, then $\tg(\tx) \leq \tg(p_r)$.

\item[(e)] If $s\increasesto r$, then $p_s\increasesto p_r$.
\end{itemize}
\end{lemma}

\begin{proof}
If $\tx < p_r$,  then $\ta(x) \not = r$ by definition.
But if $\talpha(\tx) > r$, then since $\alpha$
is degree one, there is a $\ty<\tx<p_r$ with $\talpha(\ty) = r$,
contradicting the definition of $p_r$, and so we have (a), and then (b)
follows immediately.
Now to prove (c), since $\ta\tg(p_r) = 2 \ta(p_r) = 2 r$, again by the
definition of $p_r$, we have $\tg(p_r) \geq p_{2r}$. Now if
$x \leq p_r$ and $\tg(\tx) > p_{2r}$, there would be a $\ty < x \leq p_r$
with $\tg(\ty) = p_{2r}$. But then $2 \ta(\ty) = \ta\tg(\ty) = \ta(p_{2r})
= 2r$, and so $\ta(\ty) = r$, contradicting the definition of
$p_r$. Thus $x \leq p_r$ implies $\tg(x) \leq p_{2r}$,
so we have (c), and then immediately (d).
Finally, if $s\nearrow r$, by (b), $\{p_s\}$ is increasing in $s$ and
is bounded above by $p_r$. If there was a $z < p_r$ with $p_s\nearrow z$,
then by the continuity of $\ta$, $\ta(z) = r$, again contradicting the
definition of $p_r$.
\end{proof}

For $r = k/2^n$ with $k\in\integers$ and $n \in \naturals$
we adapt the special notation of $p_{k,n} = p_r$. 
By conjugating $g$ by  a rigid rotation we may assume 
that $p_{0, 0} = 0$, which since $\alpha$ is degree one implies  
$p_{k, 0} = k$ for all $k$, and so using Lemma~\ref{th5}c,
$\tg^n(\pkn) = k$ for all $k,n$.
The next lemma gives an explicit consequence of 
a non-injective semiconjugacy in the form of 
a ``fold'' in the dynamics of $g$.

\begin{lemma}\label{th6}
If  $g:\circle\raw \circle$ is a continuous,
degree-two circle map which has been conjugated so
$p_{0, 0} = 0$ and is such that its semiconjugacy $\alpha$ is light
but not injective,
then there  exists $N, K\in\naturals$ with $0 \leq K < 2^N$ and 
$\hx\in\reals$ with $\pKN < \hx < \pKNN$ so that 
$\tg^N(\hx) = K-1$.
\end{lemma}

\begin{proof} 
First note that there exists some $p_{r'}$ 
and an $x'\in\reals$ with  $x' > p_{r'}$ and $\ta(x') < r'$, for otherwise
by  Lemma~\ref{th5}ab, $\ta$ would be injective. If we fix this $x'$, then
the set $\{ r : x' > p_r \ \hbox{and}\ \ta(x') < r \}$ is nonempty.
Let $r_0$ be its supremum and note that by  Lemma~\ref{th5}e, $x' > p_{r_0}$
and $\ta(x') < r' \leq r_0$.
Next we prove that $ s > r_0$ implies $x' < p_s$, by assuming
to the contrary that $s > r_0$ and $x' \geq p_s$. Now if
$x' = p_s$, then $\ta(x') = s > r_0$,  and if
$x' > p_s$, by the definition of $r_0$ we have $\ta(x') \geq s
> r_0$. Thus in either case we have a contradiction to $\ta(x') < r_0$.

Letting $s_0 = \ta(x')$, since $s_0 < r_0$, elementary number theory   
yields integers $K$ and $N$ with
\begin{equation*}
2^N s_0 + 1 + 2\|\varphi\| < K \leq 2^N r_0 < K+1.
\end{equation*}
with $\varphi$ as in \eqref{cc}.
Then since $2^N s_0 = 2^N \ta(x') = \ta \tg^N(x')$, \eqref{gs} says that
$| \tg^N(x') - 2^N s_0 | < \|\varphi \|$ and so $\tg^N(x') < K-1$. 
Now since $K/2^N \leq r_0 < (K+1)/2^N$, using the first paragraph of
the proof and  Lemma~\ref{th5}b, we have
$p_{K,N} \leq p_{r_0} < x' < p_{K+1, N}$. By hypothesis
$p_{0, 0} = 0$, and so by  Lemma~\ref{th5}c, $\tg^N(p_{K,N}) = K$ 
and  $\tg^N(p_{K+1,N}) = K+1$, and thus
$K-1 \in \tg^N([p_{K,N},p_{K+1,N}]) $. 
The continuity of $\tg$ then yields the required $\hx$.
Since $\ta(x + 1) = \ta(x) + 1$, we may assume that 
$0 \leq K/2^N < 1$. \end{proof}

\section{The main theorem}

The main theorem gives a number of conditions which are
equivalent to $g$ having a light semiconjugacy that
is not injective.  It will easily imply Theorem~\ref{th2} of  the introduction.

\begin{theorem}\label{th7}
If  $g$ is a continuous,
 degree-two circle map with a light semiconjugacy
$\alpha$, then  the following are equivalent:
\begin{itemize}

\item[(a)] The map $\tg$ is not injective,
 
\item[(b)] The map $\alpha$ is not injective,
 
\item[(c)] There exists a full measure, dense, $G_\delta$-set
$\Lambda\subset\circle$
so that $\theta\in \Lambda$ implies that $\alpha\Inv(\theta)$ is uncountable,
and thus contains a Cantor set,

\item[(d)] The topological entropy of $g$ satisfies $\htop(g) > \log(2)$,

\item[(e)] For all nontrivial intervals $J \subset\circle$,
the map $\alpha\vert_J$ is not of bounded variation.
\end{itemize}
\end{theorem}

\begin{proof}
If $\alpha$ is injective, then so is $\tg = \ta \td \ta\Inv$,
and thus  (a) implies (b). Since conjugate maps have the
same entropy, (d) implies (b). 
Now if $g$ is injective, then by \eqref{aa}, $\ta$ is nondecreasing,
but by hypothesis   $\alpha$ is light, and so $\ta$ is strictly increasing
and thus is injective, therefore (b) implies (a).
The fact that each of (c) and (e) imply (b) is obvious, so we
henceforth assume that $\ta$ is not injective and show that
this implies (c), (d), and (e).

Let $K$ and $N$ be as in Lemma~\ref{th6} and continue to assume that
$g$ has been conjugated so that $p_{0,0} = 0$.  By  Lemma~\ref{th6} 
and  Lemma~\ref{th4}a, we may find
intervals $I_a$, $I_b$, and $I_c$ in $\Ispec$ with disjoint interiors and
$I_a \leq I_b \leq I_c$
so that $\tg^N(I_a) =\tg^N(I_b) 
= [K-1, K]$ and $\tg^N(I_c) = [K, K+1]$. For each $k \not = K$ with
$0 \leq k < 2^N$, define intervals
$I_k = [p_{k,N}, p_{k+1, N}]$.  Define a set of ``addresses'' as
$A = \{0, 1, 2, \ \dots \ , K-1, K+1, \ \dots \ 2^N-1, a, b, c\}$, and for
$\eta\in A$,  let $\phi(\eta)$ be given by $\phi(a) = \phi(b) = K-1$,
$\phi(c) = K$, and for $0 \leq k < 2^N$, $\phi(k) = k$. 
By  Lemma~\ref{th5}c we now have that for all $\eta\in A$, $\tg^N(I_\eta)$
covers $[0, 1] + \phi(\eta)$. Projecting the $I_\eta$ to the circle
we see that $g^N$ has a $(2^N + 2)$-fold horseshoe and
so (see Theorem 4.3.2 in \cite{ALM}),  $\htop(g^N) \geq \log(2^N + 2)$ and
therefore  $\htop(g) \geq \log(2^N + 2)/N > \log(2)$,  yielding (d).

Returning to the covering space $\reals$, since $g$ is a degree-two map,
for any integer $m$, $\tg^N(I_\eta + m)$ covers $[0, 1] + \phi(\eta) + 2^N m$.
Thus by  Lemma~\ref{th4}b 
for any sequence $\us \in A^\naturals$ we may find a $y\in [0,1]$ with 
\begin{equation}\label{sum}
\tg^{Nj}(y) \in I_{s_j} + \sum_{i=0}^{j-1} 2^{N(j -i -1)} \phi(s_i)
\end{equation}
for all $j\in\naturals$. Now a given $y$ can represent two or
more sequences,
but that can only happen if for some $i$, $\tg^{Ni}(y)$ is contained
in two intervals and so must be in the boundary of some $I_\eta$, but
then by construction of the $I_\eta$, $\tg^{N(i+1)}(y) \in \integers$, 
and since $p_{0,0} = 0$ as noted above  Lemma~\ref{th6} we have
that for all $j> i$,   $\tg^{Nj}(y) \in \integers$. If we assume
initially that that $K \not = 0, 2^N - 1$, then for any integer $m$,
$(I_{2^N -1} + m -1) \cap (I_0 + m)= \{m\}  $. Thus a point $y$ 
can represent two sequences $\us$ and $\us'$ only if 
$s_j$ and $s_j^\prime$ are contained
in $\{2^N-1, 0\}$ for all sufficiently large $j$.
Therefore, if we say a sequence has
a \de{nontrivial tail} if there exist arbitrarily large $j$
with $s_j\not\in \{2^N-1, 0\}$, we see that when $\us$ has
a nontrivial tail, $\us \not = \us'$ implies that the
corresponding $y's$ are distinct. To make this true
when $K=0$
the definition of nontrivial tail must be altered to 
to require arbitrarily large $j$ with $s_j\not\in \{2^N-1, a\}$, 
and when $K=2^N-1$
to require arbitrarily large $j$ with $s_j\not\in \{c, 0\}$. 

Now note that \eqref{aa} implies that a $y$ which satisfies \eqref{sum}
will have 
\begin{equation*}
\ta(y) = \lim_{j\raw\infty} 
\frac{1}{2^{Nj}}\sum_{i=0}^{j-1} 2^{N(j -i -1)} \phi(s_i)
= \sum_{i=0}^{\infty} \frac{\phi(s_i)}{2^{N(i+1)}}.
\end{equation*}
Since $\phi(a) = \phi(b) = \phi(K-1) = K-1$, whenever
$\phi(s_i) = K-1$ in this sum there are three possible choices of
$s_i$ which give the same value of $\ta(y)$. Thus if 
 $\ut\in\{ 1, 2, \ \dots \ 2^N -1 \}^\naturals$ is a 
sequence with $t_i  = K-1$ for infinitely many $i$,  the sum
\begin{equation}\label{rsum}
r = \sum_{i=0}^{\infty} \frac{t_i}{2^{N(i+1)}} 
\end{equation}
is equal to the sum in \eqref{sum} for uncountably many sequences
$\us$. If uncountably many of these sequences $\us$ have
nontrivial tail, then for such an $r$ the set $\ta\Inv(r)$ is uncountable.
We will prove that the collection of all such $r$ is
as in (c).  

It is well known that when a map is ergodic
with respect to a smooth measure on a compact manifold, the collection
of $x$ whose orbits are dense is a  dense, $G_\delta$, full measure
set, and that the angle-doubling map
$d$ is ergodic with respect to Lebesgue measure.  Thus
(c) is proven after we show that whenever $\theta$ has a dense
orbit, its lift to an $r\in [0,1)$ is as described at the
end of the previous paragraph.

The proof of this proceeds by repeating the construction that gave rise to 
\eqref{sum}  in the easier
case of $\td$. For $0\leq k < 2^N$, let $\hI_k = [k/2^N, (k+1)/2^N]$, and
so for any integer $m$, 
$\td^N(\hI_k + m)$ covers $[0, 1] + k + 2^N m$. Thus 
for any sequence
$\ut\in\{ 1, 2, \ \dots \ 2^N -1 \}^\naturals$ 
we may find an $r\in [0,1)$
with 
\begin{equation}\label{tsum}
\td^{Nj}(r) \in \hI_{t_j} + \sum_{i=0}^{j-1} 2^{N(j -i -1)} t_i 
\end{equation}
for all $j\in\naturals$. This implies that $r$ is given by \eqref{rsum}. 
Conversely, because $\td$ is expanding and for all $k,m$, 
$\td^N(I_k + m)=[0, 1] + k + 2^N m$, it follows that  
any $r\in[0,1)$ with $d^{Nj}(r) \not\in\integers$ for all $j$ will be the
 unique $r\in[0,1]$ satisfying \eqref{tsum} for a sequence 
$\ut$ with $t_i\not\in \{0, 2^N-1\}$ for arbitrarily large i.
In particular, if $\theta\in\circle$ has a dense orbit under $d$, then its
orbit lands infinitely often in the projection to the
circle of every interval $\hI_k$, and thus its lift
$r\in[0,1)$ yields a sequence $\ut$ for which $t_i = K-1$ infinitely often, 
and any $\us$ with $\phi(t_i) = s_i$ for all $i$ must have
a nontrivial tail. Thus for such $r$,
$\ta\Inv(r)$ is uncountable and 
thus $\alpha\Inv(\theta)$ is uncountable also, proving (c).

Now to prove (e), say that an interval $J\subset \reals$ \de{unit covers}
if for some integer $M$, $[M, M+1] \subset J$. By construction for each
$\eta\in A$, $\tg^N(I_\eta)$ unit covers. Since there are $2^N + 2$
such intervals $I_\eta$, using $\ta = (\ta\circ\tg^N)/2^N$ we obtain
that the variation of $\ta$ on the interval $[0,1]$ satisfies
$\var(\alpha, [0, 1]) \geq  (2^N + 2)/2^N$. Now since each  $\tg^N(I_\eta)$
unit covers and each unit interval $[M, M+1]$ contains $I_\eta + M$
for all $\eta$, using  Lemma~\ref{th4} there are $2^N +2$ intervals $I_{\eta,j}$
in each $I_\eta$ so that each $\tg^{2N}(I_{\eta,j})$ unit covers, 
so $\var(\alpha, [0, 1]) \geq  (2^N + 2)^2/2^{2N}$. An obvious induction
then yields that for all $j$, 
\begin{equation}\label{varest}
\var(\alpha, [0, 1]) \geq  \frac{(2^N + 2)^j}{2^{Nj}},
\end{equation}
which goes to infinity as $j\raw\infty$, and
so $\ta$ has unbounded variation on $[0,1]$. 

Now as noted at the beginning of the proof of  Lemma~\ref{th3}, 
for any interval $J\subset\reals$ there 
is a $w\in \naturals$ so that $\tg^w(J)$ unit covers. Then using
 Lemma~\ref{th4} and the intervals of the previous paragraph we get that
\begin{equation*}
\var(\ta, J) \geq  \frac{(2^N + 2)^j}{2^{Nj + w}} \raw \infty,
\end{equation*}
proving (e).\end{proof}
\medskip
\noindent\textit{Proof of Theorem~\ref{th2}.}
Assume that the  semiconjugacy of $g$ to $d$  has monotone-light
decomposition, $\alpha = \ell m$. Now if $J$ is an interval
such that $m(J)$  is a point, then certainly $d \ell m(J) =  \ell  m g(J)$ is
also a point, and since $\ell$ is light, this says 
that $mg(J)$ must also be a point. 
Thus the formula $\hg = m g m\Inv$ unambiguously defines
a continuous degree-two map with light conjugacy $\ell$.
Now if $\alpha$ has a disconnected preimage, then 
$\ell$ must also, and so by Theorem~\ref{th7}c, (a) implies (b),
while the converse is trivial.
The graph of $\ell$ differs from that of $\alpha$ only by
the insertion of perhaps a countable number of horizontal intervals, and
so assuming (a), by Theorem~\ref{th7}c, (c) follows and the converse
is also clear. Finally, since $g$ is semiconjugate to $\hg$,
$\htop(g) \geq \htop(\hg)$.
Assuming (a), Theorem~\ref{th7}d gives $\htop(\hg) > \log(2)$, finishing
the proof. \hfill $\Box$

\section{Remarks and questions}

\begin{remark} The primary distinction between the general case of
Theorem~\ref{th2} and the light semiconjugacy case of Theorem~\ref{th7}
is that a general $g$ can have have an arbitrary
amount of dynamical complications and thus entropy in, 
say, a periodic interval. In this case, one can have
$\htop(g) > \log(2)$,
which clearly implies that $\alpha$ is not injective, but
it doesn't necessarily imply that $\alpha$ is not monotone.
\end{remark}

\begin{remark} If $g$ is piecewise monotone with a finite number
of turning points, 
then it follows from a theorem of Misiurewicz and Szlenk (\cite{MS}) that
the variation estimate on $\tg^{Nj}$ that gives rise to 
\eqref{varest} is equivalent to the entropy result.
\end{remark}

\begin{remark}\label{r3} As noted in  the introduction, there is a theorem in
symbolic dynamics which has similarities with Theorems~\ref{th2} 
and ~\ref{th7}. This
theorem says that if $(\Sigma, \sigma)$ and $(\Sigma', \sigma')$ 
are transitive  subshifts of finite type with $\htop(\sigma) > 0$, 
and $\alpha$ is a surjective semiconjugacy, then there is a dichotomy:
either  $\htop(\sigma) = \htop(\sigma')$ and there exists an
integer $N$ so that the cardinality of every $\alpha\Inv(\us)$ is 
at most $N$, or $\htop(\sigma) > \htop(\sigma')$ and
for the topologically generic point $\us'\in \Sigma'$, the 
point inverse $\alpha\Inv(\us')$ is uncountable.  
See Corollary 4.1.8 in \cite{K} as well as Section 6 in \cite{A}.

The proof of Theorem~\ref{th7} given here has much of 
the flavor of this symbolic dynamics result, basically showing
the existence of a diamond in the semiconjugacy. In fact,
parts of the result could have been reduced to the symbolic dynamics theorem,
but doing so would have resulted in a longer, less self-contained
proof.
\end{remark}

\begin{remark} Parts of Theorem~\ref{th7} can also be obtained by more
topological methods. From Proposition~\ref{th3} it follows that
any $g$ with a light conjugacy is locally eventually onto,
and from this is follows fairly easily that if $\alpha$ is
not injective in one open set, then it is not injective
in any open set. Such an $\alpha$ is called nowhere
locally injective. Block, Oversteegen and Tymchatyn have shown that
any light, nowhere locally injective map
between manifolds has the property that
the topologically generic point has 
a Cantor set as its point inverse (\cite{BOT}).
\end{remark}

\begin{remark} This paper has dealt primarily with combinatorial/topological
aspects of degree-two circle maps. It would also be of
interest to study quantitative/analytic aspects. For example,
for a $g$ with a light semiconjugacy,
give an explicit relationship between properties
of its semiconjugacy $\alpha$, say the fractal dimensions of 
the graph of $\alpha$, and the difference in entropy, $\htop(g) -
\log(2)$. In this regard we note then when $g$ has
a finite number of turning points, its semiconjugacy
can be treated in the context of fractal functions.
In particular, if $g$ is piecewise linear
with expanding pieces, then $\alpha$ is an affine fractal function
and its graph is the attractor of a planar iterated
function system,   (see \cite{Pu}). 
Also, in analogy to the degree-one case,
it would also be interesting to study the
transition to a non-monotone semiconjugacy
in parameterized families, for example in the standard
degree-two family $f_{b,\omega}(x) = 2x + \omega + b\sin(2\pi x)$.
\end{remark}

\providecommand{\bysame}{\leavevmode\hbox to3em{\hrulefill}\thinspace}
\providecommand{\MR}{\relax\ifhmode\unskip\space\fi MR }
\providecommand{\MRhref}[2]{%
  \href{http://www.ams.org/mathscinet-getitem?mr=#1}{#2}
}
\providecommand{\href}[2]{#2}


\begin{thebibliography}{10}

\bibitem{A}
Roy~L. Adler, \emph{Symbolic dynamics and {M}arkov partitions}, Bull. Amer.
  Math. Soc. (N.S.) \textbf{35} (1998), no.~1, 1--56.

\bibitem{ALM}
Llu{\'{\i}}s Alsed{\`a}, Jaume Llibre, and Micha{\l} Misiurewicz,
  \emph{Combinatorial dynamics and entropy in dimension one}, Advanced Series
  in Nonlinear Dynamics, vol.~5, World Scientific Publishing Co. Inc., River
  Edge, NJ, 2000.

\bibitem{BC}
L.~S. Block and W.~A. Coppel, \emph{Dynamics in one dimension}, Lecture Notes
  in Mathematics, vol. 1513, Springer-Verlag, Berlin, 1992.

\bibitem{BOT}
A.~Blokh, L.~Oversteegen, and E.~Tymchatyn, \emph{On almost one-to-one maps},
  Trans. Amer. Math. Soc., to appear.

\bibitem{dMvS}
Welington de~Melo and Sebastian van Strien, \emph{One-dimensional dynamics},
  Ergebnisse der Mathematik und ihrer Grenzgebiete (3) [Results in Mathematics
  and Related Areas (3)], vol.~25, Springer-Verlag, Berlin, 1993.

\bibitem{F}
John Franks, \emph{Anosov diffeomorphisms}, Global Analysis (Proc. Sympos. Pure
  Math., Vol. XIV, Berkeley, Calif., 1968), Amer. Math. Soc., Providence, R.I.,
  1970, pp.~61--93.

\bibitem{H}
Michael Handel, \emph{Entropy and semi-conjugacy in dimension two}, Ergodic
  Theory Dynam. Systems \textbf{8} (1988), no.~4, 585--596.

\bibitem{K}
Bruce~P. Kitchens, \emph{Symbolic dynamics}, Universitext, Springer-Verlag,
  Berlin, 1998.

\bibitem{Pu}
Peter~R. Massopust, \emph{Fractal functions, fractal surfaces, and wavelets},
  Academic Press Inc., San Diego, CA, 1994.

\bibitem{MS}
M.~Misiurewicz and W.~Szlenk, \emph{Entropy of piecewise monotone mappings},
  Studia Math. \textbf{67} (1980), no.~1, 45--63.

\bibitem{MN}
Micha{\l} Misiurewicz and Zbigniew Nitecki, \emph{Combinatorial patterns for
  maps of the interval}, Mem. Amer. Math. Soc. \textbf{94} (1991), no.~456,
  vi+112.

\end{thebibliography}
\end{document}